\documentclass{article}
\usepackage{amssymb,latexsym,amsmath,amsthm,mathrsfs}
\usepackage{graphicx}
\newcommand{\comment}[1]{}

\begin{document}
\title{A general method for investigating the roots of all equations
by approximation\footnote{Presented to the St. Petersburg Academy on
April 25, 1776. Originally published as
{\em Methodus generalis investigandi radices omnium aequationum per approximationem},
Nova Acta Academiae Scientarum Imperialis Petropolitinae \textbf{6} (1790),
16--24. E643 in the Enestr{\"o}m index.
Translated from the Latin by Jordan Bell, School of Mathematics and
Statistics, Carleton University, Ottawa, Canada.
Email: jbell3@connect.carleton.ca}}
\author{Leonhard Euler}
\date{}
\maketitle

\S 1. If $Z$ is an unknown quantity whose value is taken from an equation, it
can always be exhibited under the form $Z=0$, and thus $Z$ is a certain
function of $z$.
Therefore a determinate value is searched for, which if written in place of
$z$
turns this function $Z$ into nothing.
Just as if this equation were given: $z^3=3z+2$,
placing all the terms on the same side it will be $z^3-3z-2=0$,
and thus the given function $Z=z^3-3z-2$. Then a value is searched for,
which when substituted
for $z$ makes this function $Z$ equal to nothing, as clearly happens in
this case by taking $z=2$.

\S 2. If, however, in such a proposed equation $Z=0$,
any other value $v$ is written in place of $z$, from which the function
$Z$ takes a value $=V$, then at least it will not be $V=0$: for
if it came out $V=0$, this symbol $v$ would be a true root $z$.
Thus in the announced example, $Z=z^3-3z-2$,
if in place of $z$ we write $3$, so that it will become $v=3$, it would be
$V=16$, and in no way $=0$.

\S 3. Therefore since in the given equation $Z=0$, if one writes
in place of $z$ any value $v$, which makes $Z=V$, and it is not $V=0$,
let us put $V=y$, where $y$ will be some known function of $v$,
which vanishes if for $v$ the true value of the root $z$
were taken.

\S 4. Since therefore $y=V$ is some function of $v$,
it comprises a certain curve, where at the abscissa $=v$ corresponds
the ordinate
$=y$. Therefore permuting the coordinates, so that now $y$ exhibits the
abscissa and $v$ the ordinate, this ordinate $v$ can likewise be seen
as a certain function of $y$, and it will therfore be the case
that if one takes $y=0$ then the ordinate $v$ will turn out to be equal
to the root $z$ that is sought, and so the entire question is reduced to this:
what value does this function of $y$, namely $v$, take if one writes
$0$ in place of $y$? For then this value of $v$ will yield a true root
$z$.

\S 5. Because, therefore, we consider $v$ as a function of $y$,
we put in the general way $v=\Gamma :y$; we also
note
 from the nature of differentials that
\[
\Gamma :(y+a)=\Gamma :y+\frac{a\partial \cdot \Gamma :y}{\partial y}+
\frac{aa\partial^2\cdot \Gamma :y}{1\cdot 2\cdot \partial y^2}+
\frac{a^3\partial^3 \cdot \Gamma :y}{1\cdot 2\cdot 3\cdot \partial y^3}+
\textrm{etc.}
\]
and so by writing $v$ in place of $\Gamma :y$ we will have
\[
\Gamma :(y+a)=v+
\frac{a\partial v}{\partial y}+
\frac{aa\partial \partial v}{1\cdot 2\cdot \partial y^2}+
\frac{a^3 \partial^3 v}{1\cdot 2\cdot 3\cdot \partial y^3}+
\frac{a^4\partial^4 v}{1\cdot 2\cdot 3\cdot 4\cdot \partial y^4}+\textrm{etc.},
\]
in which expression the element $\partial y$ is taken as constant.
Or, since $v$ is a function of $y$, if we set
$\frac{\partial v}{\partial y}=p$; $\frac{\partial p}{\partial y}=q$;
$\frac{\partial q}{\partial y}=r$; $\frac{\partial r}{\partial y}=s$;
and thus for the higher to infinity,
we will obtain for finite quantities
\[
\Gamma :(y+a)=v+ap+\frac{1}{2}aaq+\frac{1}{6}a^3r+\frac{1}{24}a^4s+\frac{1}{120}a^5t+\textrm{etc.}
\]

\S 6. Because here any quantity can be assumed in place of $a$,
we may assume $a=-y$, and then $\Gamma :(y+a)$ will exhibit to us the same
function which would result from the form $\Gamma :y$ if $y+a$ were
written in place of $y$, that is $y-y=0$. We have found that $v$ will
become the root $z$ that is sought, so it would thus be $z=\Gamma :(y-y)$:
therefore, if in the series that has been found we write $-y$ in place of
$a$, we will obtain
\[
z=v-py+\frac{1}{2}qyy-\frac{1}{6}ry^3+\frac{1}{24}sy^4-\frac{1}{120}ty^5+
\textrm{etc.}
\]

\S 7. So we can make this series more accomodating for use, we 
shall take the quantity $y$ out of it. For since
$y=V$ and $V$ is a known function of $v$,
from its character we have at once $p=\frac{\partial v}{\partial V}$,
and with this value found it will be in turn $q=\frac{\partial p}{\partial V}$,
and then further $r=\frac{\partial q}{\partial V}; s=\frac{\partial r}{\partial V}$, and thus in turn to infinity. With this observed, we will be obtain
the following infinite series for the root that is sought:
\[
z=v-pV+\frac{1}{2}qV^2-\frac{1}{6}rV^3+\frac{1}{24}sV^4-\frac{1}{120}tV^5+\textrm{etc.}
\]
It will be worthwhile to illustrate the use of this series with several examples.

\begin{center}
{\Large Example I.}
\end{center}

\S 8. Let this quadratic equation be given: $zz=a$, and
one thus searches for a series which is equal to $z=\surd a$.
Here it will therefore be $Z=zz-a$, and then by writing $v$ in place of $z$ 
it will be $V=vv-a$, and hence $\partial V=2v\partial v$.
\[
p=\frac{1}{2v};\quad q=-\frac{1}{4v^3};\quad r=\frac{3}{8v^5};\quad
s=-\frac{3\cdot 5}{16v^7}; \quad t=\frac{3\cdot 5\cdot 7}{32v^9}; \quad \textrm{etc.}
\]
and having found these, we will obtain for the root that is sought
\[
\begin{split}
\surd a=&v-\frac{(vv-a)}{2v}-\frac{(vv-a)^2}{2\cdot 4\cdot v^3}-
\frac{3(vv-a)^3}{6\cdot 8\cdot v^5}\\
&-\frac{3\cdot 5\cdot (vv-a)^4}{24\cdot 16\cdot v^7}-
\frac{3\cdot 5\cdot 7\cdot (vv-a)^5}{120\cdot 32\cdot v^9}-\textrm{etc.}
\end{split}
\]

But let us therefore take this form for the proposed equation:
$zz=bb+c$,
thus so that $z=\surd(bb+c)$,
one should write everywhere $bb+c$ in place of $a$.
Then indeed because the quantity $v$ is at our discretion,
we may set $v=b$,
whence it will be $vv-a=-c$,
and with this done the root that is sought will be
\[
\surd(bb+c)=b+\frac{c}{2b}-\frac{cc}{2\cdot 4\cdot b^5}+\frac{3c^3}{6\cdot 8\cdot b^5}
-\frac{3\cdot 5\cdot c^4}{24\cdot 16\cdot b^7}+\frac{3\cdot 5\cdot 7\cdot c^5}{120\cdot 32\cdot b^9}-\textrm{etc.}.
\]
This
series is obtained by the extraction of the root by the binomial expansion.
And thus if $\surd 10$ is sought,
one can take $b=3$ \& $c=1$, from which gathers 
\[
\surd 10=3+\frac{1}{6}-\frac{1}{226}+\frac{1}{3888}-\textrm{etc.}
\]
Therefore it will be approximately $\surd 10=3+\frac{1}{6}$,
and the series will be made more highly convergent by putting
$b=3\frac{1}{6}=\frac{19}{6}$;
then also it will be $c=-\frac{1}{36}$,
whence by only the first two terms one will find
\[
\surd 10=3\frac{1}{6}-\frac{1}{12\cdot 19}=3\frac{37}{228},
\]
whose square comes  so close to $10$ that the error is only
$\frac{1}{51984}$.

\begin{center}
{\Large Example II.}
\end{center}

\S 9. Let this equation be proposed: $z^n=a$,
thus so that $z=\sqrt[n]{a}$ is searched for.
It will therefore be $Z=z^n-a$, and so $V=v^n-a$ and
$\partial V=nv^{n-1}\partial v$;
whereby the letters $p,q,r,s$, etc. will be determined in the following way:
\[
\begin{split}
&p=\frac{1}{nv^{n-1}};\quad q=-\frac{(n-1)}{nnv^{2n-1}}; \quad r=\frac{(n-1)(2n-1)}{n^3 v^{3n-1}};\\
&s=-\frac{(n-1)(2n-1)(3n-1)}{n^4 v^{4n-1}}; \quad \textrm{etc.}
\end{split}
\]
and having substituted these values, the root that is sought will be
\[
\begin{split}
&z=\sqrt[n]{a}=v-\frac{(v^n-a)}{nv^{n-1}}-\frac{(n-1)(v^n-a)^2}{2nnv^{2n-1}}
-\frac{(n-1)(2n-1)(v^n-a)^3}{6n^3v^{3n-1}}\\
&-\frac{(n-1)(2n-1)(3n-1)(v^n-a)^4}{24n^4v^{4n-1}}-\textrm{etc.}
\end{split}
\]
Namely this series will always express the same value, whatever number is taken
in place of $v$, which of course remains entirely at our discretion.
However, this series will approach the truth more quickly the smaller
the difference between the power $z^n$ and the number $a$.

But if it had been $a=b^n+c$, it would be convenient to take
$v=b$, for then it would be $v^n-a=-c$ and the series expressing the root
that is sought would be
\[
\begin{split}
&\sqrt[n]{b^n+c}=b+\frac{c}{nb^{n-1}}-\frac{(n-1)cc}{2nnb^{2n-1}}+
\frac{(n-1)(2n-1)c^3}{6n^3b^{3n-1}}\\
&-\frac{(n-1)(2n-1)(3n-1)c^4}{24n^4b^{4n-1}}+
\frac{(n-1)(2n-1)(3n-1)(4n-1)c^5}{120n^5b^{5n-1}} \, \textrm{etc.}
\end{split}
\]
This series is also found from the expansion of the binomial $(b^n+c)^{\frac{1}{n}}$.

\begin{center}
{\Large Example III.}
\end{center}

\S 10. With the equation proposed: $z^3=z-1$,
to find a series, which exhibits the value
of the root $z$.
Here it will therefore be $Z=z^3-z+1$, and so
$V=v^3-v+1$,
from which it will be $\partial V=\partial v(3vv-1)$;
from this we will have
\begin{eqnarray*}
p&=&\frac{1}{3vv-1};\\
q&=&-\frac{6v}{(3vv-1)^3};\\
r&=&\frac{6(15vv+1)}{(3vv-1)^5};\\
s&=&-\frac{36v(6vv+1)}{(3vv-1)^7};\\
&&\textrm{etc.}
\end{eqnarray*}
Then substituting these values it can be gathered:
\[
\begin{split}
&z=\sqrt[3]{(z-1)}=v-\frac{(v^3-v+1)}{3vv-1}-\frac{3v(v^3-v+1)^2}{(3vv-1)^3}\\
&-\frac{(15vv+1)(v^3-v+1)^3}{(3vv-1)^5}-\frac{15v(6vv+1)(v^3-v+1)^4}{(3vv-1)^7}
-\textrm{etc.}
\end{split}
\]
where again the quantity $v$ can be assumed at our pleasure;
it will be convenient however to take it such that it does not differ
greatly from the value of the root $z$.
Thus if we take $v=1$, it follows $v^3-v+1=1$ and $3vv-1=2$,
from which it would be $z=1-\frac{1}{2}-\frac{3}{8}-\frac{1}{8}$,
which 

\S 11. From this example it is clear moreover, that 
this approximation will hardly succeed unless a value close
to the root $z$,
when the usual method does the work much better.
Nonetheless though, the series that have been found by this method are
very remarkable, because continued to infinity they always exhibit the same
value $z$, whatever number is assumed for $v$.
Also, this method especially merits all attention,
because not only does the series exhibit a root, but also
any power at all, which will be dealt with in the following
problem.

\begin{center}
{\Large Problem.}
\end{center}

{\em Given any equation $Z=0$, where $Z$ is any function of $z$,
to find an infinite series by which not only the root $z$, but also
each power $z^n$ of it can be expressed.}

\begin{center}
{\Large Solution.}
\end{center}

\S 12. One starts the calculation as before, namely, in place of $z$
is written any quantity $v$,
where $Z$ turns into $V$; and as much as $v$ is not a root of the
proposed equation, so $V$ will not turn into nothing. We then
put, as above, $V=y$, \& moreover $v$ will be seen as a function of $y$.
It thus happens that in the case in which one puts $y=0$, it will be
$v=z$,
since in this case $v$ will be a true root of the equation. 

\S 13. Indeed in a similar way also any power $v^n$ can be seen as a function
of $y$,
which therefore, when $y=0$, or when in place of $y$ is written $y-y$,
bestows to us the value that is sought $z^n$. Thus, as we have shown above,
it will be
\[
z^n=v^n-\frac{y\partial \cdot v^n}{\partial y}
+\frac{yy\partial^2\cdot v^n}{2\partial v^2}
-\frac{y^3\partial^3\cdot v^n}{2\cdot 3\cdot \partial y^3}+
\frac{y^4 \partial^4\cdot v^n}{2\cdot 3\cdot 4\cdot \partial y^4}-\textrm{etc.}
\]
Then because $y=V$, it will be
\[
z^n=v^n-\frac{V\partial \cdot v^n}{\partial V}+\frac{V^2\partial^2 \cdot v^n}{2\partial V}-\frac{V^3\partial^3\cdot v^n}{2\cdot 3\cdot \partial V^3}+
\frac{V^4\partial^4 \cdot v^n}{2\cdot 3\cdot 4\cdot \partial V^4}-\textrm{etc.}
\]

\S 14. So that we can free this series from differentials here,
let us put $P=\frac{\partial \cdot v^n}{\partial V}=\frac{nv^{n-1}\partial v}{\partial V}$;
hence, like before, the following values are formed:
\[
Q=\frac{\partial P}{\partial V}; \quad R=\frac{\partial Q}{\partial V};
\quad \frac{\partial R}{\partial V};\quad T=\frac{\partial S}{\partial V};
\quad \textrm{etc.}
\]
with these values found, the series sought that is sought expressing $z^n$ will be
\[
z^n=v^n-VP+\frac{1}{2}V^2Q-\frac{1}{6}V^3R+\frac{1}{24}V^4S-\frac{1}{120}V^5T+\textrm{etc.}
\]
where it should be noted that after all these terms are expanded in $v$,
then clearly
in place of $v$ any quantity can be assumed and still the same
value of the power $z^n$ will always result from it.

\S 15. Since this series should express the value of $z^n$, putting $n=0$
our series should produce unity,
which is clear taking the first term $v^n$ by itself;
then because the letter $P$ has a factor $n$, moreover all the following
letters will be multiplied by this same exponent $n$, and thus putting
$n=0$ all the following terms of this series will immediately vanish.

\S 16. Here we observe that besides all the powers of $z$, we can also
exhibit
its hyperbolic logarithm, namely $lz$, 
by infinite series of this type. For since all the terms, after the first,
of the series that have been found have the factor $n$,
we may write in place of them $n\Omega$, so that it would be $z^n=v^n+n\Omega$,
and in turn $\Omega=\frac{z^n-v^n}{n}$. The differential of this equation
will be
\[
\partial \Omega=z^{n-1}\partial z-v^{n-1}\partial v,
\]
and thus in the case $n=0$ it will be $\partial \Omega=\frac{\partial z}{z}-
\frac{\partial v}{v}$. The integration of this equation yields
$\Omega=l\frac{z}{w}$,
from which we gather $lz=lv+\Omega$,
where in all the terms that form $\Omega$ it is put $n=0$.

\S 17. But if therefore $\Omega$ denotes the sum of all the terms in
the case $n=0$,
as we have assumed, 
and
if $e$ denotes that number whose hyperbolic logarithm is $=1$, it will
be $z=ve^\Omega$, therefore $z=v^n e^{n\Omega}$.
Then also, if the exponential formula $e^{n\Omega}$ is converted
into an infinite series in the usual way, it will become
\[
z^n=v^n(1+n\Omega+\frac{1}{2}nn\Omega^2+\frac{1}{6}n^3\Omega^3+\frac{1}{24}n^4\Omega^4+\textrm{etc.})
\]
This should therefore be equal to the former series, which was found earlier.
Since it is necessary that expanding the letter $\Omega$
yields the series found before, it will be helpful to show it by example.

\begin{center}
{\Large Example.}
\end{center}

\S 18. With the proposed equation $z^\lambda=a$, to find a series
both for any other power $z^m$ and for its hyperbolic logarithm. Here
we will therefore have $Z=z^\lambda-a$,
and so $V=v^\lambda-a$,
from which $\partial V=\lambda v^{\lambda-1}\partial v$. Hence now the
letters introduced before, $P,Q,R$, etc. will be determined in the following way:
\begin{eqnarray*}
P&=&\frac{n}{\lambda}v^{n-\lambda};\\
Q&=&\frac{n}{\lambda}\cdot \frac{n-\lambda}{\lambda}\cdot v^{n-2\lambda};\\
R&=&\frac{n}{\lambda}\cdot \frac{n-\lambda}{\lambda}\cdot \frac{n-2\lambda}{\lambda}\cdot v^{n-3\lambda};\\
S&=&\frac{n}{\lambda}\cdot \frac{n-\lambda}{\lambda}\cdot \frac{n-2\lambda}{\lambda}\cdot \frac{n-3\lambda}{\lambda}\cdot v^{n-4\lambda}; \textrm{etc.}
\end{eqnarray*} 
Having substituted these values, the desired series will be
\[
\begin{split}
z^n=&v^n-\frac{n}{\lambda}\cdot v^{n-\lambda}(v^\lambda-a)+\frac{n}{\lambda}\cdot \frac{(n-\lambda)}{2\lambda}\cdot v^{n-2\lambda}(v^\lambda-a)^2\\
&-\frac{n}{\lambda}\cdot \frac{(n-\lambda)}{2\lambda}\cdot \frac{(n-2\lambda)}{3\lambda} v^{n-3\lambda}(v^\lambda-a)^3\\
&+\frac{n}{\lambda}\cdot \frac{n-\lambda}{2\lambda}\cdot \frac{n-2\lambda}{3\lambda}\cdot \frac{n-3\lambda}{4\lambda}\cdot v^{n-4\lambda}(v^\lambda -a)^4-\textrm{etc.}
\end{split}
\]
By putting here $n=0$ it will be
\[
\begin{split}
lz=&lv-\frac{1}{\lambda}v^{-\lambda}(v^\lambda-a)-\frac{1}{2\lambda}v^{-2\lambda}(v^\lambda-a)^2\\
&-\frac{1}{3\lambda}v^{-3\lambda}(v^\lambda-a)^3
-\frac{1}{4\lambda}v^{-4\lambda}(v^\lambda-a)^4-\textrm{etc.}
\end{split}
\]

\end{document}